\RequirePackage{float}
\documentclass[pdflatex,sn-basic]{sn-jnl}
\usepackage{graphicx}
\usepackage{multirow}
\usepackage{amsmath,amssymb,amsfonts}
\usepackage{amsthm}
\usepackage{mathrsfs}
\usepackage{placeins}
\usepackage{xcolor}
\usepackage{textcomp}
\usepackage{manyfoot}
\usepackage{booktabs}
\usepackage{array,tabularx}
\usepackage{algorithm}
\usepackage{algorithmicx}
\usepackage{algpseudocode}
\usepackage{listings}
\usepackage{enumitem}
\usepackage{microtype}
\setlist[enumerate]{font=\normalfont}
\AtBeginDocument{}
\theoremstyle{thmstyleone}
\newtheorem{theorem}{Theorem}
\newtheorem{proposition}[theorem]{Proposition}
\newtheorem{corollary}{Corollary}
\newtheorem{lemma}{Lemma}
\theoremstyle{thmstyletwo}

\theoremstyle{thmstylethree}

\newtheorem{assumption}{Assumption}
\DeclareMathOperator*{\argmin}{arg\,min}

\newcommand{\compactdisplays}{%
  \setlength{\abovedisplayskip}{5pt plus2pt minus1pt}%
  \setlength{\belowdisplayskip}{5pt plus2pt minus1pt}%
  \setlength{\abovedisplayshortskip}{1pt plus1pt}%
  \setlength{\belowdisplayshortskip}{3pt plus1pt minus1pt}}
\AddToHook{cmd/normalsize/after}{\compactdisplays}
\AddToHook{cmd/small/after}{\compactdisplays}
\AtBeginDocument{\compactdisplays}
\setlist[enumerate]{topsep=3pt,itemsep=0pt,parsep=0pt,partopsep=0pt}
\setlist[itemize]{topsep=3pt,itemsep=0pt,parsep=0pt,partopsep=0pt}
\makeatletter
\g@addto@macro\th@thmstyleone{\thm@preskip=7pt plus2pt minus1pt\thm@postskip=7pt plus2pt minus1pt}
\g@addto@macro\th@thmstyletwo{\thm@preskip=7pt plus2pt minus1pt\thm@postskip=7pt plus2pt minus1pt}
\g@addto@macro\th@thmstylethree{\thm@preskip=7pt plus2pt minus1pt\thm@postskip=7pt plus2pt minus1pt}
\AtBeginDocument{%
  \renewenvironment{proof}[1][\proofname]{%
    \par\removelastskip\pushQED{\qed}%
    \small\normalfont\topsep4pt plus1pt minus1pt\relax\trivlist
    \item[\hskip\labelsep\itshape #1\@addpunct{}]\ignorespaces
  }{\popQED\endtrivlist\@endpefalse}}
\makeatother

\begin{document}
\title[Auto-conditioned Gradient Descent Methods]{Auto-conditioned Gradient Descent Methods for Multiobjective Optimization}

\author[1]{\fnm{Shuaichenji} \sur{Jiang}}
\affil[1]{\orgdiv{National Center for Applied Mathematics in Chongqing},
  \orgname{Chongqing Normal University},
  \orgaddress{\city{Chongqing}, \country{People's Republic of China}}}

\abstract{This paper proposes two linesearch-free algorithms based on steepest descent methods for smooth multiobjective optimization: auto-conditioned steepest descent method for multiobjective optimization (AC-SDM) and auto-conditioned normalized gradient descent method for multiobjective optimization (AC-NGDM). Both algorithms estimate the local Lipschitz constant from the function values and the first-order information obtained in previous iterations. AC-SDM uses a shared estimate of the Lipschitz constant to construct a scaled steepest descent direction. For AC-NGDM, we use each function's local Lipschitz estimate to normalize the gradient of the corresponding function and construct the descent direction. We prove that every accumulation point of any sequence generated by both algorithms is a Pareto critical point. We also establish convergence rates of $O(1/\sqrt{k})$ for nonconvex objectives and $O(1/k)$ for convex objectives. Numerical results on several test problems demonstrate that AC-SDM and AC-NGDM require fewer iterations and less CPU time than the multiobjective steepest descent method with  an Armijo line search and a max-type nonmonotone line search.}


\maketitle

\section{Introduction}
	
	In this paper, we consider the unconstrained multiobjective optimization
	problem
	\begin{equation}\label{min}
		\min_{x\in\mathbb{R}^n}
		f(x):=\bigl(f_1(x),\ldots,f_m(x)\bigr)^T,
	\end{equation}
	where each objective function $f_i:\mathbb{R}^n\to\mathbb{R}$ is continuously
	differentiable and $L_i$-smooth. 
	
	Multiobjective optimization is widely used in many research areas, which can be
	found in
 machine learning~\citep{bischof2025loss,chen2024tradeoff}, multi-energy system
	design~\citep{terlouw2023energy}, and electrocatalyst design~\citep{xu2024discovering}.
	
Over the past two decades, multiobjective first-order methods have been widely studied,  which include steepest descent methods~\citep{fliege2000}, projected gradient methods~\citep{drummond2004projected}, proximal gradient methods~\citep{tanabe2019proximal}, and nonlinear conjugate gradient methods~\citep{lucambio2018}. Although these methods construct descent directions in different ways, they commonly use line search to determine a suitable stepsize.

Line search is commonly used the backtracking strategy to select a stepsize satisfying a prescribed decrease condition. This backtracking process may require multiple evaluations of all objective functions, thereby increasing the computational cost. When applied to multiobjective optimization, line search requires a stepsize to satisfy the decrease conditions for all objective functions, which may lead to a very small stepsize~\citep{mita2019,chen2025improvements}. \citet{mita2019} proposed nonmonotone line search strategies to alleviate the small-stepsize issue, although repeated backtracking still be required. To avoid the backtracking process, researchers have recently begun exploring adaptive stepsize strategies. \citet{hoai2025adaptive} proposed four adaptive stepsize rules based on the objective gradient-difference. \citet{goncalves2026adaptive} developed a adaptive stepsize rule based on accumulated squared norms of multiobjective steepest descent directions.

We recently noted that \citet{lan2024projected} proposed an auto-conditioned projected gradient method that uses both function values and gradient information to estimate the local Lipschitz constant for stepsize selection. They also introduced a strategy to address potential nonmonotonicity in function values. Then \citet{yagishita2025simple} extended this approach to proximal gradient methods for nonconvex optimization. Inspired by these works, we adopt their approach to estimating local Lipschitz constants. We propose two algorithms, AC-SDM uses a shared estimate to scale the steepest descent direction and AC-NGDM uses each objective estimates to normalize the gradients in the direction-finding subproblem. This normalization also enables AC-NGDM to handle the imbalance caused by large differences in gradient norms across objectives. Furthermore, we establish its global convergence and convergence rate.
	
The main contributions of this paper are summarized as follows.
\begin{itemize}
\item We propose AC-SDM and AC-NGDM, which require neither line search nor prior knowledge of the global Lipschitz constants. By local Lipschitz estimates, we establish sufficient descent properties and convergence for both methods. We also establish an $O(k^{-1/2})$ convergence rate for nonconvex problems and an $O(k^{-1})$ convergence rate for convex problems for both algorithms.
\item Numerical experiments compare both methods with the multiobjective steepest descent method using an Armijo line search (MSD--Armijo) and a max-type nonmonotone line search (MSD-max type).

\end{itemize}

The paper is organized as follows. Section~2 introduces the notation and optimality concepts. Sections~3 and~4 present AC-SDM and AC-NGDM, respectively, and establish their convergence properties and convergence rates. Section~5 reports numerical comparisons with MSD--Armijo and MSD-max type. Section~6 concludes the paper and discusses future work.

	\section{Preliminaries}
	\subsection{Notations and definitions}

We use $\langle \cdot,\cdot \rangle$ and $\|\cdot\|$ to represent the standard inner product and Euclidean norm in $\mathbb{R}^n$, respectively. For any positive integer $m$, let $[m]:=\{1,\ldots,m\}$. The relations $\le$ and $<$ in $\mathbb{R}^m$ are understood componentwise; that is, $u\le v$ (respectively, $u<v$) if and only if $u_i\le v_i$ (respectively, $u_i<v_i$) for all $i\in[m]$.

The simplex in $\mathbb{R}^m$ is given by $\Delta_m=\left\{\lambda\in\mathbb{R}^m:
\sum_{i\in[m]}\lambda_i=1,\ 
\lambda_i\ge 0,\ i\in[m]\right\}.$
	We denote the Jacobian matrix and the gradient of $f_i$ at $x$ by ${J}f(x) \in \mathbb{R}^{m\times n}$ and $\nabla f_i(x) \in \mathbb{R}^n$, respectively. Then each objective function
	$f_i:\mathbb{R}^n\to\mathbb{R}$ is assumed to be $L_i$-smooth, meaning that its gradient is Lipschitz continuous with constant $L_i$, namely,
	\[
	\|\nabla f_i(x)-\nabla f_i(y)\|
	\le L_i\|x-y\|,
	\qquad
	\forall x,y\in\mathbb{R}^n,\quad i\in[m].
	\]
	By the descent lemma, this condition implies that
\begin{equation}\label{eq:descent-lemma-intro}
      f_i(y)\leq f_i(x) +\langle\nabla f_i(x),y-x\rangle +\frac{L_i}{2}\|y-x\|^2,\quad \forall x,y\in\mathbb{R}^n,\quad i\in[m].
    \end{equation}
	Let $L:=\max_{i\in[m]}L_i.$
	Since $L_i\le L$ for every $i\in[m]$, each objective function is also $L$-smooth.
	We now recall the notions of optimality for problem (\ref{min}). A point $x^*\in\mathbb{R}^n$ is called a weak Pareto solution if there is no $x\in\mathbb{R}^n$ such that $f(x)<f(x^*)$. It is called a Pareto solution if there is no $x\in\mathbb{R}^n$ such that $f(x)\le f(x^*)$ and $f(x)\neq f(x^*)$. We also say that $\bar{x}\in\mathbb{R}^n$ is \textit{Pareto stationary} (or \textit{critical}) if
	$$
	\max_{i\in[m]}\langle \nabla f_i(\bar{x}),d\rangle\ge 0,\qquad \forall d\in\mathbb{R}^n.
	$$
	
	\subsection{The Steepest Descent Method for Multiobjective Optimization}
	We first review the steepest descent method for multiobjective optimization, as presented in~\citep{fliege2000}.
	For a given $x \in \mathbb{R}^n$, we consider the following subproblem (\ref{subproblem}):
	\begin{equation} \label{subproblem}
		\min\limits_{d \in \mathbb{R}^n} \max\limits_{i\in[m]} \langle \nabla f_i(x),d\rangle+\frac{1}{2} \|d\|^2.
	\end{equation}
	Since $\langle \nabla f_i(x),d\rangle+\frac{1}{2}\|d\|^2$ is strongly convex for each $i \in [m]$, problem (\ref{subproblem}) admits a unique minimizer. We denote by $\theta(x)$ and $v(x)$ the optimal value and the optimal solution of (\ref{subproblem}), respectively. Therefore,
	\begin{equation} \label{theta}
		\theta(x)=\min\limits_{d \in \mathbb{R}^n} \max\limits_{i\in[m]} \langle \nabla f_i(x),d\rangle+\frac{1}{2} \|d\|^2,
	\end{equation}
	and
	\begin{equation} \label{v}
		v(x)=\argmin\limits_{d \in \mathbb{R}^n} \max\limits_{i\in[m]} \langle \nabla f_i(x),d\rangle+\frac{1}{2} \|d\|^2.
	\end{equation}
	An equivalent formulation with improved solvability can be obtained from the original problem:
	\begin{equation*} \tag{$QP$} \label{QP}
		\min\limits_{(t,d)\in \mathbb{R} \times \mathbb{R}^n} t+\frac{1}{2} \|d\|^2;\qquad \text{s.t.} \langle \nabla f_i(x),d \rangle \leq t, \qquad i \in [m]
	\end{equation*}
	Since problem (\ref{QP}) is convex and possesses a feasible point $(t,d)=(0,0)$, we have $\theta(x) \leq 0$. By Slater's condition, strong duality holds; hence the dual problem has the same optimal value as (QP)~\citep{rockafellar1970convex}. We consider the dual problem:
	$$\max\limits_{\lambda \in \Delta_m}\min\limits_{(t,d)\in\mathbb{R}\times\mathbb{R}^n}L(\lambda,(t,d)),$$
	where $L(\lambda,(t,d))$ is the Lagrangian,
	$$L(\lambda,(t,d))=t+\frac{1}{2}\|d\|^2+\sum\limits_{i\in[m]}\lambda_i(\langle \nabla f_i(x),d\rangle-t).$$
	From the KKT conditions~\citep{rockafellar1970convex}, we obtain
	\begin{equation} \label{kkt}
		v(x)=-\sum\limits_{i\in[m]}\lambda_i(x)\nabla f_i(x),
	\end{equation}
	where $\lambda(x)=(\lambda_1(x),...,\lambda_m(x))$ is an optimal solution of the following problem (\ref{dual}).
	\begin{equation} \label{dual}
		-\min\limits_{\lambda\in\Delta_m}\frac{1}{2}\|\sum\limits_{i\in [m]}\lambda_i\nabla f_i(x)\|^2.
	\end{equation}
	By strong duality, we have:
	$$\theta(x)=-\frac{1}{2}\|\sum\limits_{i\in [m]}\lambda_i\nabla f_i(x)\|^2=-\frac{1}{2}\|v(x)\|^2\leq0.$$
	If the $\lambda_{i}(x) \neq 0$ then
	$$\langle \nabla f_i(x),v(x)\rangle=-\|v(x)\|^2.$$
	
	Then, we will present some properties of $\theta(x)$ and $v(x)$.
	\begin{lemma}[{\citep[Lemma~3.2]{tanabe2019proximal}}] \label{thetax}
		Let $\theta(x)$ and $v(x)$ be as in (\ref{theta}) and (\ref{v}),
		respectively. Then,
		\begin{enumerate}[font=\upshape]
			\item[(a)] the following conditions are equivalent:
			\begin{enumerate}[font=\upshape]
				\item[(i)] The point $x$ is Pareto critical;
				\item[(ii)] $\theta(x)=0$;
				\item[(iii)] $v(x)=0$.
			\end{enumerate}
			\item[(b)] the following conditions are equivalent:
			\begin{enumerate}[font=\upshape]
				\item[(i)] The point $x$ is non-critical;
				\item[(ii)] $\theta(x)<0$;
				\item[(iii)] $v(x)\neq0$.
			\end{enumerate}
			\item[(c)] the functions $v(x)$ and $\theta(x)$ are continuous for all $x \in \mathbb{R}^n$.
		\end{enumerate}
	\end{lemma}
\begin{lemma}[{\normalfont\citep[Lemma~1]{fliege2000}}]
	\label{inequation}
	Let $\theta(x)$ and $v(x)$ be defined by \eqref{theta} and \eqref{v},
	respectively.
	\begin{enumerate}[label=(\alph*)]
		\item If $x$ is Pareto critical, then $v(x)=0\in\mathbb{R}^n$ and
		$\theta(x)=0$.
		
		\item If $x$ is not Pareto critical, then
		$$
		\theta(x)\leq-\frac{1}{2}\|v(x)\|^2<0,
		~and\quad	\langle\nabla f_i(x),v(x)\rangle\leq\theta(x),
			\quad i=1,\ldots,m.
		$$
	\end{enumerate}
\end{lemma}

	\section{Auto-conditioned steepest descent method for multiobjective optimization}
	 We apply (\ref{eq:descent-lemma-intro}) to estimate the local Lipschitz constant, and use it to scale the steepest descent direction. Then we use a constant $\alpha$, which will be used to ensure a descent property of the objective functions in the subsequent analysis.
	\begin{algorithm}[htbp]
		\caption{Auto-conditioned steepest descent method for multiobjective optimization (AC-SDM)}\label{AC_noconvex}
		\begin{algorithmic}[1]
			\Require Initial point $x^0\in\mathbb{R}^n$, parameters $\alpha>1$, $L_0>0$ and $\varepsilon>0$.
			\State $\gamma_1=L_0$.
			\For{$k=1,2,\ldots$}
			\State
			$$
			\lambda^k\in\argmin_{\lambda\in\Delta_m}
			\left\|\sum_{i=1}^{m}\lambda_i\nabla f_i(x^{k-1})\right\|,
            \qquad v^k=-\sum_{i=1}^{m}\lambda_i^k\nabla f_i(x^{k-1}).
			$$
			\If{$\|v^k\|\leq\varepsilon$}
			\State \Return $x^{k-1}$.
			\EndIf
			\State
			$$
			d^k=\frac{1}{\alpha\gamma_k}v^k,
			\qquad
			x^k=x^{k-1}+d^k.
			$$
			\State For all $i\in[m]$,
			\begin{equation} \label{the definion of Li}
L_i^k=
\frac{2\bigl(f_i(x^k)-f_i(x^{k-1})
  -\langle\nabla f_i(x^{k-1}),d^k\rangle\bigr)}
  {\|d^k\|^2}.
			\end{equation}
			\State
			$$
			L^k=\max_{i\in[m]}L^k_i,
			\qquad
			\gamma_{k+1}=\max\{\gamma_k,L^k\}.
			$$
			\EndFor
		\end{algorithmic}
	\end{algorithm}
	
	The auto-conditioned direction in Algorithm~1 is the unique solution of the following subproblem, defined for $x\in\mathbb{R}^n$, $\gamma>0$, and $\alpha>1$:
	\begin{equation} \label{newsub}
		\min\limits_{d\in \mathbb{R}^n}\max\limits_{i \in [m]}\langle \nabla f_i(x),d \rangle+\frac{\alpha\gamma}{2}\|d\|^2.
	\end{equation}
	Thus, $d^k$ is a positive rescaling of the steepest descent direction $v^k$. Next we establish the relationship between $v(x)$ and $d(x)$. 
	\begin{corollary}
		Let $x\in\mathbb{R}^n$, $\gamma>0$, and $\alpha>1$,
		and let $d(x)$ be defined by
		$$
		d(x)
		=
		\argmin_{d\in\mathbb{R}^n}
		\max_{i\in[m]}
		\langle\nabla f_i(x),d\rangle
		+
		\frac{\alpha\gamma}{2}\|d\|^2.
		$$
		Then the following statements are equivalent:
		\begin{enumerate}
			\item[(i)] $d(x)=0$;
			\item[(ii)] $v(x)=0$;
			\item[(iii)] $x$ is Pareto critical.
		\end{enumerate}
	\end{corollary}
	\begin{proof}
		By \eqref{kkt}, we have
		$$
		d(x)
		=
		-\frac{1}{\alpha\gamma}
		\sum_{i=1}^{m}\lambda_i\nabla f_i(x)
		=
		\frac{1}{\alpha\gamma}v(x).
		$$
		Since $\alpha\gamma>0$, statements \textnormal{(i)} and
		\textnormal{(ii)} are equivalent. By Lemma \ref{thetax}, statements
		\textnormal{(ii)} and \textnormal{(iii)} are equivalent. Then suppose that $x$ is Pareto critical. There exists $\mu\in\Delta_m$ such that
	$
		\sum_{i=1}^m \mu_i\nabla f_i(x)=0.
	$
		Define
	$
		\widehat{\lambda}_i
		:=\frac{\gamma_i\mu_i}{\sum_{j=1}^m\gamma_j\mu_j},
		~ i\in[m].
	$
		Since $\gamma_i>0$, we have $\widehat{\lambda}\in\Delta_m$ and
		\[
		\sum_{i=1}^m\frac{\widehat{\lambda}_i}{\gamma_i}\nabla f_i(x)
		=
		\frac{\sum_{i=1}^m\mu_i\nabla f_i(x)}
		{\sum_{j=1}^m\gamma_j\mu_j}
		=0.
		\]
		Therefore, the three statements are equivalent.
	\end{proof}

	\begin{theorem}\label{thm:acsdm-finite}
		Let $\{\gamma_k\}$ and $\{L^k\}$ be generated by Algorithm~\ref{AC_noconvex}.
		Let $\mathcal{S}=\{k\geq 1:\beta \gamma_k\geq L^k\}$, where $\beta=\frac{\alpha+1}{2}>1$. Then there are only finitely many indices $k\notin \mathcal{S}$. Consequently, $L^k\leq\beta\gamma_k$ for all sufficiently large $k$.
	\end{theorem}
	\begin{proof}
		It is clear that $L^k_i\leq L$ for all $i\in [m]$ and $k\geq 1$. Hence,
		$$L_0\leq\gamma_k\leq\max\{L_0,L\}.$$
		If $k\notin \mathcal{S}$, then $L^k>\beta\gamma_k$, and thus $\gamma_{k+1}=L^k$. Therefore, the number of indices $k\notin \mathcal{S}$ is bounded by
		$$\left\lceil \log_{\beta} \frac{\max\{L_0, L\}}{L_0} \right\rceil. $$
	\end{proof}
	
	Combining Theorem~\ref{thm:acsdm-finite} with the descent property of the steepest descent direction, we establish the following sufficient decrease property of AC-SDM.
	
	\subsection{Global convergence}
	\begin{lemma}\label{lem:acsdm-descent}
		Let $\{x^k\}$ be the sequence produced by Algorithm~\ref{AC_noconvex}. There exist a finite integer $N_1$ and a constant $c=\frac{\alpha L_0}{2}
		\left(1-\frac{\beta}{\alpha}\right)$ such that, for all $k > N_1$, the following inequality holds:
		$$
		f_i(x^{k-1})-f_i(x^k)
		\geq
		c\|d^k\|^2,
		\quad i\in[m],\quad k>N_1.
		$$
	\end{lemma}
	\begin{proof}
		From the definition of $L_i^k$, we have
		\begin{equation} \label{fk-fk-1}
			f_i(x^k)-f_i(x^{k-1})
			=
			\frac{1}{\alpha\gamma_k}\left\langle \nabla f_i(x^{k-1}),\, v^k \right\rangle
			+
			\frac{L_i^k}{2\alpha^2\gamma_k^2}\|v^k\|^2.
		\end{equation}
		By Lemma \ref{inequation},
\[
    f_i(x^k)-f_i(x^{k-1})\leq -\frac{1}{2\alpha\gamma_k}\|v^k\|^2 +\frac{L_i^k}{2\alpha^2\gamma_k^2}\|v^k\|^2 \leq -\frac{1}{2\alpha\gamma_k} \left(1-\frac{L^k}{\alpha\gamma_k}\right)\|v^k\|^2.
    \]
		By Theorem~\ref{thm:acsdm-finite}, for any $k > N_1$ and any $i \in [m]$, it holds that $L_i^k \leq L^k \leq \beta \gamma_k$, then 
\begin{equation}\label{eq:sufficient-descent}
      f_i(x^{k-1})-f_i(x^k)\geq \frac{1}{2\alpha\gamma_k} \left(1-\frac{\beta}{\alpha}\right)\|v^k\|^2 \geq \frac{\alpha L_0}{2} \left(1-\frac{\beta}{\alpha}\right)\|d^k\|^2 \geq c\|d^k\|^2.
    \end{equation}
		
	\end{proof}
	\begin{assumption}\label{ass:bounded-level-set-ac}
	Let $\{x^k\}$ be the sequence generated by Algorithm~1, and let
	$N_1$ be the finite integer in Lemma~\ref{lem:acsdm-descent}. The lower level set
	$$
	\Omega_{N_1}
	:=
	\left\{
	x\in\mathbb{R}^n
	\mid
	f_i(x)\leq f_i(x^{N_1}),\quad i\in[m]
	\right\}
	$$
	is bounded.
\end{assumption}
	\begin{theorem}
		Let $\{x^k\}_k$ be the sequence generated by Algorithm~1.
		Suppose that Assumption~\ref{ass:bounded-level-set-ac} holds. Then
		$\{x^k\}$ is bounded, and every accumulation point of $\{x^k\}$ is a
		Pareto critical point.
	\end{theorem}
	\begin{proof}
		By Lemma~\ref{lem:acsdm-descent}, for every $i\in[m]$, the sequence
		$\{f_i(x^k)\}_{k\geq N_1}$ is nonincreasing. Moreover, for all $k\geq N_1$,
		$x^k\in\Omega_{N_1}$. Hence, Assumption~\ref{ass:bounded-level-set-ac}
		implies that $\{x^k\}_k$ is bounded and that
		$\{f_i(x^k)\}_k$ is bounded below. Therefore,
		$
		f_i(x^{k-1})-f_i(x^k)\to0.
		$
		Together with \eqref{eq:sufficient-descent}, this yields
		$
		\|v^k\|\to0.
		$
		
		Let $\{x^{k_j}\}$ be any convergent subsequence such that
		$x^{k_j}\to x^\ast$. Since
		$
		x^{k_j}-x^{k_j-1}
		=
		\frac{1}{\alpha\gamma_{k_j}}v^{k_j}
		$
		and $\gamma_{k_j}\geq L_0$, we have
		$$
		\|x^{k_j}-x^{k_j-1}\|
		\leq
		\frac{1}{\alpha L_0}\|v^{k_j}\|
		\to0.
		$$
		Thus, $x^{k_j-1}\to x^\ast$. Recalling that
		$v^{k_j}=v(x^{k_j-1})$, the continuity of $v$ gives
		$
		v(x^\ast)
		=
		\lim_{j\to\infty}v(x^{k_j-1})
		=
		0.
		$
		By Lemma~1, $x^\ast$ is a Pareto critical point.
	\end{proof}

	We now establish convergence rates for AC-SDM, considering the nonconvex and convex cases.
	\subsection{The non-convex case}
	\begin{theorem}\label{thm:ACSDM_d_rate}
		Let $\{x^k\}$ be the sequence generated by Algorithm~1, and suppose that
		Assumption~\ref{ass:bounded-level-set-ac} holds. For each $i\in[m]$, let	$f_{i,\Omega_{N_1}}^{\min} := \min_{x\in\Omega_{N_1}}f_i(x)$,
	
		Then
	$$
		\min_{N_1+1\leq \ell\leq k}\|d^\ell\|
		\leq
		\frac{
			\sqrt{\min\limits_{i\in[m]}
				\left(f_i(x^{N_1})-f_{i,\Omega_{N_1}}^{\min}\right)}
		}{
			\sqrt{c(k-N_1)}
		},
		\qquad \forall k>N_1.
	$$
	\end{theorem}
	\begin{proof}
		By Lemma~\ref{lem:acsdm-descent}, for every $k>N_1$ and every $i\in[m]$, 
		and summing the inequality from $\ell=N_1+1$ to $k$, we get
		$$
		\sum_{\ell=N_1+1}^{k}
		\bigl(f_i(x^{\ell-1})-f_i(x^\ell)\bigr)
		\ge
		c\sum_{\ell=N_1+1}^{k}\|d^\ell\|^2.
		$$
	Then
		$$
		f_i(x^{N_1})-f_i(x^k)
		\ge
		c\sum_{\ell=N_1+1}^{k}\|d^\ell\|^2.
		$$
		Since $x^k\in\Omega_{N_1}$, we have
		$$
		c\sum_{\ell=N_1+1}^{k}\|d^\ell\|^2
		\leq
		f_i(x^{N_1})-f_i(x^k)
		\leq
		f_i(x^{N_1})-f_{i,\Omega_{N_1}}^{\min}.
		$$
		Therefore,
		$$
		\min_{N_1+1\leq \ell\leq k}\|d^\ell\|
		\leq
		\frac{
			\sqrt{\min\limits_{i\in[m]}
				\left(f_i(x^{N_1})-f_{i,\Omega_{N_1}}^{\min}\right)}
		}{
			\sqrt{c(k-N_1)}
		},
		\qquad \forall k>N_1.
		$$
	\end{proof}

	\subsection{The convex case}
	In this subsection, we consider the convergence rate of Algorithm~1 under the
	convexity assumption. 

	For the subsequent convergence-rate analysis, we recall the following merit function {\citep{tanabe2024merit}}:
	$$
	u_0(x):=
	\sup_{y\in\mathbb{R}^n}
	\min_{i\in[m]}
	\{f_i(x)-f_i(y)\}.
	$$
	\begin{proposition}[{\citep[Theorem~3.1]{tanabe2024merit}}]
		$x^\ast \in \mathbb{R}^n$ is a weak Pareto solution of problem~\eqref{min} if and only if
		$u_0(x^\ast)=0$.
	\end{proposition}
	
	\begin{lemma} \label{decrease}
		Assume that all objective functions $f_i(x)$, $i \in [m]$ are convex .
		Let $\{x^k\}$ be the sequence generated by Algorithm~1, and let
		$N_1$ be the finite integer in Lemma~\ref{lem:acsdm-descent}.
		Then, for every $y\in\mathbb{R}^n$ and every $k>N_1$,
		the following inequality holds:
		$$
		\sum_{i=1}^m \lambda_i^k
		\bigl(
		f_i(x^k)-f_i(y)
		\bigr)
		\leq
		\frac{\alpha\gamma_k}{2}
		\left(
		\|x^{k-1}-y\|^2-\|x^k-y\|^2
		\right).
		$$
	\end{lemma}
	\begin{proof}
		For all $k> N_1$,
		$
		L_i^k\leq \beta\gamma_k<\alpha\gamma_k.
		$
		By (\ref{fk-fk-1}) and convexity of $f_i$, we obtain
		$$
		f_i(x^k)-f_i(y)
		\leq
		\left\langle
		\nabla f_i(x^{k-1}),
		x^{k-1}-y+d^k
		\right\rangle
		+
		\frac{\alpha\gamma_k}{2}\|d^k\|^2.
		$$
		Multiply both sides by $\lambda_i^k$ and sum over $i=1,\ldots,m$,
\[
    \sum_{i=1}^m\lambda_i^k\bigl(f_i(x^k)-f_i(y)\bigr) \leq\left\langle \sum_{i=1}^m\lambda_i^k\nabla f_i(x^{k-1}), x^{k-1}-y+d^k\right\rangle +\frac{\alpha\gamma_k}{2}\|d^k\|^2.
    \]
		Based on the definition of $d(x)$, we have
\[
    \begin{aligned}
      \sum_{i=1}^m\lambda_i^k\bigl(f_i(x^k)-f_i(y)\bigr) &\leq-\alpha\gamma_k \langle d^k,x^{k-1}-y+d^k\rangle +\frac{\alpha\gamma_k}{2}\|d^k\|^2\\
      &=-\alpha\gamma_k\langle d^k,x^{k-1}-y\rangle -\frac{\alpha\gamma_k}{2}\|d^k\|^2.
    \end{aligned}
    \]
		Using the identity
		$$
		\|(x^{k-1}-y)+d^k\|^2
		=
		\|x^{k-1}-y\|^2
		+
		2\langle d^k,x^{k-1}-y\rangle
		+
		\|d^k\|^2 ,
		$$
		the following inequality holds
		$$
		\sum_{i=1}^m \lambda_i^k
		\bigl(
		f_i(x^k)-f_i(y)
		\bigr)
		\leq
		\frac{\alpha\gamma_k}{2}
		\left(
		\|x^{k-1}-y\|^2-\|x^k-y\|^2
		\right).
		$$
	\end{proof}

	\begin{theorem}
		Assume that each $f_i$ is convex, for $i\in [m]$. Let $\{x^k\}$ be generated by Algorithm~1, and suppose that Assumption~\ref{ass:bounded-level-set-ac} holds. Let $\Theta$ be an upper bound on the diameter of $\Omega_{N_1}$. Then, for all $k>N_1$,
		$$
		u_0(x^k)
		\leq
		\frac{\alpha\max\{L_0,L\}\Theta^2}{2(k-N_1)}.
		$$
	\end{theorem}
	\begin{proof}
		Fix $k>N_1$ and define
		$$
		\Omega_k
		:=
		\{y\in\mathbb{R}^n\mid f(y)\leq f(x^k)\}.
		$$
		For any $y\in\Omega_k$ and $l=N_1,\ldots,k-1$, by Lemma~\ref{decrease} and the monotonicity of
		$\{f(x^l)\}$ we can get
\[
    \sum_{i=1}^m\lambda_i^{l+1} \bigl(f_i(x^{l+1})-f_i(y)\bigr) \leq\frac{\alpha\gamma_{l+1}}{2} \left(\|x^l-y\|^2-\|x^{l+1}-y\|^2\right).
    \]
		The left-hand side is nonnegative. Using $\gamma_{l+1}\leq\max\{L_0,L\}$ and
		$f_i(x^k)\leq f_i(x^{l+1})$. Summing over $l=N_1,\ldots,k-1$ yields
\[
    \begin{aligned}
      \sum_{l=N_1}^{k-1}\sum_{i=1}^m\lambda_i^{l+1} \bigl(f_i(x^k)-f_i(y)\bigr) &\leq\frac{\alpha\max\{L_0,L\}}{2} \left(\|x^{N_1}-y\|^2-\|x^k-y\|^2\right)\\
      &\leq\frac{\alpha\max\{L_0,L\}\Theta^2}{2},
    \end{aligned}
    \]

		On the other hand, since $\lambda^{l+1}\in\Delta_m$. Consequently,
		$$
		\min_{i\in[m]}
		\{f_i(x^k)-f_i(y)\}
		\leq
		\frac{\alpha\max\{L_0,L\}\Theta^2}{2(k-N_1)},
		\qquad y\in\Omega_k.
		$$
		Taking the supremum over $y\in\Omega_k$ and using
		$
		u_0(x^k)
		=
		\sup_{y\in\Omega_k}
		\min_{i\in[m]}
		\{f_i(x^k)-f_i(y)\},
		$
	so, we have
		$$
	u_0(x^k)
	\leq
	\frac{\alpha\max\{L_0,L\}\Theta^2}{2(k-N_1)}.
	$$
	\end{proof}

	\section{Auto-conditioned normalized gradient descent method for multiobjective optimization}
Since $\gamma_k$ retains the largest estimate across all objectives, this may lead to a very small scaling factor, thereby slowing convergence. To address this limitation, we propose a normalized variant of AC-SDM that assigns a parameter $\gamma_i^k$ to each objective and uses these to normalize the objective gradients in the direction-finding subproblem. The resulting method has the same convergence properties and rates as AC-SDM.
	
	For any $x\in\mathbb{R}^n$, different from the subproblem used in
	Algorithm~1, we consider the following scaled subproblem:
	$$
	d_\gamma(x)
	=
	\argmin_{d\in\mathbb{R}^n}
	\max_{i\in[m]}
	\left\langle
	\frac{1}{\gamma_i}\nabla f_i(x),d
	\right\rangle
	+
	\frac{\alpha}{2}\|d\|^2 .
	$$
	Here each $\gamma_i>0$ is a positive constant. By (\ref{kkt}), we have
\[
    d_\gamma(x)=-\frac{1}{\alpha}\sum_{i=1}^m \frac{\lambda_i^\gamma(x)}{\gamma_i}\nabla f_i(x),\qquad \lambda^\gamma(x)\in\argmin_{\lambda\in\Delta_m} \left\|\sum_{i=1}^m\frac{\lambda_i}{\gamma_i} \nabla f_i(x)\right\|.
    \]
	Next, we establish the relationship between $d_{\gamma}(x)$ and Pareto critical points and present several properties of $d_{\gamma}(x)$.

	\begin{theorem}\label{thm:dgamma-critical}
		Let $x\in\mathbb{R}^n$ and let $d_\gamma(x)$ be the solution
		of the scaled subproblem defined above. Then $d_\gamma(x)=0$
		if and only if $x$ is a Pareto critical point.
	\end{theorem}
	\begin{proof}
		If $d_\gamma(x)=0$, then
		$$\frac{1}{\alpha}
		\sum_{i=1}^m
		\frac{\lambda_i^\gamma(x)}{\gamma_i}
		\nabla f_i(x)
		=
		0.
		$$
		Define
		$
		\mu_i(x)
		:=
		\frac{\lambda_i^\gamma(x)/\gamma_i}
		{\sum_{j=1}^m \lambda_j^\gamma(x)/\gamma_j},
		~ i\in[m].
		$
		Since $\lambda^\gamma(x)\in\Delta_m$ and $\gamma_i>0$, we have
		$
		\mu(x)\in\Delta_m.
		$
		Moreover,
		$$
		\sum_{i=1}^m
		\mu_i(x)\nabla f_i(x)
		=
		\frac{
			\sum_{i=1}^m
			\frac{\lambda_i^\gamma(x)}{\gamma_i}
			\nabla f_i(x)}
		{\sum_{j=1}^m \lambda_j^\gamma(x)/\gamma_j}
		=
		0.
		$$
		Hence
		$
		0\in
		\operatorname{conv}
		\{\nabla f_1(x),\ldots,\nabla f_m(x)\}.
		$
		Thus $x$ is Pareto critical point.
	\end{proof}

	\begin{lemma}\label{lem:dgamma-inequation}
		Let $x\in\mathbb{R}^n$, and let $\theta_\gamma(x)$ and
		$d_\gamma(x)$ be defined by
\[
    \begin{aligned}
      \theta_\gamma(x)&=\min_{d\in\mathbb{R}^n}\max_{i\in[m]}
        \left\langle\frac{1}{\gamma_i}\nabla f_i(x),d\right\rangle
        +\frac{\alpha}{2}\|d\|^2,\\
      d_\gamma(x)&=\argmin_{d\in\mathbb{R}^n}\max_{i\in[m]}
        \left\langle\frac{1}{\gamma_i}\nabla f_i(x),d\right\rangle
        +\frac{\alpha}{2}\|d\|^2,
    \end{aligned}
    \]
		respectively, where $\gamma_i>0$ for all $i\in[m]$.
		Then the following statements hold:
		\begin{enumerate}[label=(\roman*)]
			\item If $x$ is Pareto critical, then
			$$
			d_\gamma(x)=0\in\mathbb{R}^n
			\quad\text{and}\quad
			\theta_\gamma(x)=0.
			$$
			
			\item If $x$ is not Pareto critical, then
\[
    \theta_\gamma(x)\leq-\frac{\alpha}{2}\|d_\gamma(x)\|^2<0,\qquad \left\langle\frac{1}{\gamma_i}\nabla f_i(x),d_\gamma(x) \right\rangle\leq\theta_\gamma(x),\quad i=1,\dots,m.
    \]
			\item The mapping
			$
			(x,\gamma)\mapsto d_\gamma(x)
			$
			is continuous on $\mathbb{R}^n\times\mathbb{R}_{++}^m$.
		\end{enumerate}
	\end{lemma}
	\begin{proof}
		Let
		$$
		g_i^\gamma(x):=\frac{1}{\gamma_i}\nabla f_i(x),
		\qquad i\in[m].
		$$
	The scaled subproblem has the same form as \eqref{subproblem}, with
		$\nabla f_i(x)$ replaced by $g_i^\gamma(x)$. Hence, statement ~\textnormal{(i)}~\textnormal{(ii)} follow directly from Lemma~\ref{inequation}.
		Statement~\textnormal{(iii)} follows from the definition of the problem defining
		$d_\gamma(x)$.
	\end{proof}
	
	\begin{algorithm}[H]
		\caption{Auto-conditioned normalized gradient descent method for multiobjective optimization (AC-NGDM)} \label{ac_convex}
		\begin{algorithmic}[1]
			\Require Initial point $x^0\in\mathbb{R}^n$, parameters $\alpha>1$, $L_0>0$ and $\varepsilon>0$.
			\State $\gamma_i^1=L_0$, $i\in[m]$.
			\For{$k=1,2,\ldots$}
			\State
\[
      d_\gamma^k=-\frac{1}{\alpha}\sum_{i=1}^m \frac{(\lambda_\gamma^k)_i}{\gamma_i^k} \nabla f_i(x^{k-1}),\qquad \lambda_\gamma^k\in\argmin_{\lambda\in\Delta_m} \left\|\sum_{i=1}^m\frac{\lambda_i}{\gamma_i^k} \nabla f_i(x^{k-1})\right\|.
      \]
		
			\If{$\|d_\gamma^k\|\leq\varepsilon$}
			\State \Return $x^{k-1}$.
			\EndIf
			\State
			$$
			x^k=x^{k-1}+d_\gamma^k.
			$$
			\State For all $i\in[m]$,
			\begin{equation}
L_i^k=
\frac{2\bigl(f_i(x^k)-f_i(x^{k-1})
  -\langle\nabla f_i(x^{k-1}),d_\gamma^k\rangle\bigr)}
  {\|d_\gamma^k\|^2}.
			\end{equation}
			\State
			$$
			\gamma_i^{k+1}=\max\{\gamma_i^k,L_i^k\},
			\qquad i\in[m].
			$$
			\EndFor
		\end{algorithmic}
	\end{algorithm}
	\begin{theorem}\label{thm:acngdm-finite}
		Let $\{\gamma_i^k\}$ and $\{L_i^k\}$, $i\in[m]$, be generated by Algorithm~\ref{ac_convex}.
		Let $\mathcal{S}_i=\{k\geq 1:\beta \gamma_i^k\geq L_i^k\}$ and $\mathcal{S}_\gamma=\cap_{i=1}^m \mathcal{S}_i$, where $\beta=\frac{\alpha+1}{2}>1$. Then there are only finitely many indices $k\notin \mathcal{S}_\gamma$. Consequently, $L_i^k\leq\beta\gamma_i^k$ for all $i\in[m]$ and all sufficiently large $k$.
	\end{theorem}
	\begin{proof}
		It is clear that $L^k_i\leq L$ for all $i\in [m]$ and $k\geq 1$. Hence,
		$$L_0\leq\gamma_i^k\leq\max\{L_0,L_i\}.$$
		If $k\notin \mathcal{S}_i$, then $L_i^k>\beta\gamma_i^k$, and thus $\gamma_i^{k+1}=L_i^k$. Therefore, the number of indices $k\notin \mathcal{S}_i$ is bounded by
		$$\left\lceil \log_{\beta} \frac{\max\{L_0, L_i\}}{L_0} \right\rceil. $$
		Thus, since $m$ is finite, there are only finitely many indices $k\notin \mathcal{S}_\gamma$.
	\end{proof}
	
AC-NGDM also satisfies an sufficient decrease property.
	\begin{lemma}\label{lem:acngdm-descent}
		Let $\{x^k\}$ be the sequence produced by Algorithm~2. Then there exist a
		finite integer $N_2$ and a constant
		$c_\gamma
		=
		\frac{\alpha L_0}{2}
		\left(
		1-\frac{\beta}{\alpha}
		\right)>0$,
		such that, for all $k>N_2$ and all $i\in[m]$, the following inequality holds:
		\begin{equation}\label{eq:acngdm-sufficient-descent}
			f_i(x^{k-1})-f_i(x^k)
			\geq
			c_\gamma\|d_\gamma^k\|^2,
			\quad i\in[m],\quad k>N_2.
		\end{equation}
	\end{lemma}
	\begin{proof}
		By the definition of $L_i^k$, the update rule
		and Lemma~\ref{lem:dgamma-inequation},
		Thus,
\[
    f_i(x^k)-f_i(x^{k-1})\leq-\frac{\alpha\gamma_i^k}{2}\|d_\gamma^k\|^2 +\frac{L_i^k}{2}\|d_\gamma^k\|^2 =-\frac{\alpha\gamma_i^k}{2} \left(1-\frac{L_i^k}{\alpha\gamma_i^k}\right) \|d_\gamma^k\|^2.
    \]
		By Theorem~\ref{thm:acngdm-finite}, there exists a finite integer $N_2$ such that
\[
    \begin{aligned}
      f_i(x^{k-1})-f_i(x^k)\geq\frac{\alpha L_0}{2} \left(1-\frac{\beta}{\alpha}\right)\|d_\gamma^k\|^2 \geq c_\gamma\|d_\gamma^k\|^2, \qquad\forall k>N_2,\quad i\in[m].
    \end{aligned}
    \]
	\end{proof}

\begin{assumption}\label{ass:bounded-level-set-nac}
	Let $\{x^k\}$ be the sequence generated by Algorithm~2, and let
	$N_2$ be the finite integer in Lemma~\ref{lem:acngdm-descent}. The lower level set
	$$
	\Omega_{N_2}
	:=
	\left\{
	x\in\mathbb{R}^n
	\mid
	f_i(x)\leq f_i(x^{N_2}),\quad i\in[m]
	\right\}
	$$
	is bounded.
\end{assumption}
	\begin{lemma}\label{lem:acngdm-critical}
		Let $\{x^k\}$ be the sequence produced by Algorithm~2. Suppose that
		Assumption~\ref{ass:bounded-level-set-nac} holds. Then every accumulation point
		of $\{x^k\}$ is a Pareto critical point.
	\end{lemma}
	\begin{proof}
	By Lemma~\ref{lem:acngdm-descent}, there exist $N_2<\infty$ and
	$c_\gamma>0$ such that, for every $k>N_2$ and $i\in[m]$,
	$$
	f_i(x^{k-1})-f_i(x^k)
	\geq c_\gamma\|d_\gamma^k\|^2.
	$$
	Therefore, $x^k\in\Omega_{N_2}$ for all $k\geq N_2$. Since $\Omega_{N_2}$ is
	compact, summing the above inequality gives
	$
	\sum_{k=N_2+1}^{\infty}\|d_\gamma^k\|^2<\infty,
	$
	and hence
	$$
	d_\gamma^k\to0,
	\qquad
	x^k-x^{k-1}=d_\gamma^k\to0.
	$$
	Moreover, for each $i\in[m]$, $\{\gamma_i^k\}$ is nondecreasing and
	$
	L_0\leq\gamma_i^k\leq\max\{L_0,L_i\}.
	$
	Thus, $\gamma^k\to\gamma^\ast\in\mathbb{R}_{++}^m$.
	Let $x^\ast$ be an accumulation point and choose
	$x^{k_j}\to x^\ast$. Then $x^{k_j-1}\to x^\ast$. By
	Lemma~\ref{lem:dgamma-inequation} and the continuity of
	$(x,\gamma)\mapsto d_\gamma(x)$,
	$$
	d_{\gamma^\ast}(x^\ast)
	=
	\lim_{j\to\infty}
	d_{\gamma^{k_j}}(x^{k_j-1})
	=
	\lim_{j\to\infty}d_\gamma^{k_j}
	=0.
	$$
	Therefore, $x^\ast$ is a Pareto critical point. 
\end{proof}

	\subsection{The non-convex case}
	\begin{theorem}\label{thm:ACNGDM_d_rate}
		Let $\{x^k\}$ be the sequence generated by Algorithm~2, and suppose that
		Assumption~\ref{ass:bounded-level-set-nac} holds. For each $i\in[m]$, let
		$$
		f_{i,\Omega_{N_2}}^{\min}
		:=
		\min_{x\in\Omega_{N_2}}f_i(x),
		\quad
		M_\gamma:=\frac{\alpha L_0}{2}
		\left(1-\frac{\beta}{\alpha}\right).
		$$
		Then
		$$
		\min_{N_2+1\leq \ell\leq k}\|d_\gamma^\ell\|
		\leq
		\frac{
			\sqrt{\min\limits_{i\in[m]}
				\left(f_i(x^{N_2})-f_{i,\Omega_{N_2}}^{\min}\right)}
		}{
			\sqrt{M_\gamma(k-N_2)}
		},
		\qquad \forall k>N_2.
		$$
	\end{theorem}
	\begin{proof}
		By Lemma~\ref{lem:acngdm-descent} and
		$\gamma_i^\ell\geq L_0$, for every $i\in[m]$, we have
		$$
		M_\gamma
		\sum_{\ell=N_2+1}^{k}\|d_\gamma^\ell\|^2
		\leq
		f_i(x^{N_2})-f_i(x^k)
		\leq
		f_i(x^{N_2})-f_{i,\Omega_{N_2}}^{\min}.
		$$
		Here, the last inequality follows from $x^k\in\Omega_{N_2}$. Combining the above
		inequalities and taking the minimum over
		$i\in[m]$ yield the desired result.
			$$
		\min_{N_2+1\leq \ell\leq k}\|d_\gamma^\ell\|
		\leq
		\frac{
			\sqrt{\min\limits_{i\in[m]}
				\left(f_i(x^{N_2})-f_{i,\Omega_{N_2}}^{\min}\right)}
		}{
			\sqrt{M_\gamma(k-N_2)}
		},
		\qquad \forall k>N_2.
		$$
	\end{proof}

	\subsection{The convex case}
	
In this subsection, we establish the convergence rate of Algorithm~2 in the
convex case. The relevant properties and their proofs are similar to those for Algorithm~1.
	\begin{lemma}\label{lem:acngdm-convex-decrease}
		Assume that all objective functions $f_i$, $i\in[m]$, are convex.
		Let $\{x^k\}$ be the sequence generated by Algorithm~2, and let
		$N_2$ be the finite integer in Lemma~\ref{lem:acngdm-descent}. Then, for all
		$k>N_2$, the following inequality holds:
\[
    \sum_{i=1}^m\frac{\lambda_i^k}{\gamma_i^k} \bigl(f_i(x^k)-f_i(y)\bigr) \leq\frac{\alpha}{2} \left(\|x^{k-1}-y\|^2-\|x^k-y\|^2\right),\quad \forall y\in\mathbb{R}^n.
    \]
	\end{lemma}
	\begin{proof}
		This proof is same as Lemma~\ref{decrease}. The only change is that
		Algorithm~2 satisfies
		$$
		d_\gamma^k
		=
		-\frac{1}{\alpha}
		\sum_{i=1}^m
		\frac{\lambda_i^k}{\gamma_i^k}
		\nabla f_i(x^{k-1}),
		\qquad
		x^k=x^{k-1}+d_\gamma^k.
		$$
		Using convexity and $L_i^k\leq \beta\gamma_i^k<\alpha\gamma_i^k$ for all
		$k>N_2$, we get
		$$
		f_i(x^k)-f_i(y)
		\leq
		\left\langle
		\nabla f_i(x^{k-1}),
		x^{k-1}-y+d_\gamma^k
		\right\rangle
		+
		\frac{\alpha\gamma_i^k}{2}\|d_\gamma^k\|^2.
		$$
		Multiplying by $\lambda_i^k/\gamma_i^k$, summing over $i$, and using the
		definition of $d_\gamma^k$, we obtain
\[
    \sum_{i=1}^m\frac{\lambda_i^k}{\gamma_i^k} \bigl(f_i(x^k)-f_i(y)\bigr) \leq\frac{\alpha}{2} \left(\|x^{k-1}-y\|^2-\|x^k-y\|^2\right),\quad \forall y\in\mathbb{R}^n.
    \]
	\end{proof}

	\begin{theorem}\label{thm:acngdm-convex-rate}
		Assume that each $f_i$ is convex, for $i\in[m]$. Let $\{x^k\}$ be generated by
		Algorithm~2, and suppose that Assumption~\ref{ass:bounded-level-set-nac}
		holds. Let $\Theta$ be an upper bound on the diameter of
		$\Omega_{N_2}$. Then, for all $k>N_2$,
		the following inequality holds:
		$$
		u_0(x^k)
		\leq
		\frac{\alpha\max\{L_0,L\}\Theta^2}{2(k-N_2)}.
		$$
	\end{theorem}
	\begin{proof}
		The proof is the same as that of the convex convergence rate theorem for
		Algorithm~1, with Lemma~\ref{lem:acngdm-convex-decrease} replacing
		Lemma~\ref{decrease}. We only indicate the required modification.
		
		Fix $k>N_2$ and let
		$\Omega_k:=\{y\in\mathbb{R}^n:f(y)\leq f(x^k)\}$.
		For any $y\in\Omega_k$, summing
		Lemma~\ref{lem:acngdm-convex-decrease} from $l=N_2$ to $k-1$ gives
		$$
		\sum_{l=N_2}^{k-1}
		\sum_{i=1}^m
		\frac{\lambda_i^{l+1}}{\gamma_i^{l+1}}
		\bigl(
		f_i(x^k)-f_i(y)
		\bigr)
		\leq
		\frac{\alpha\Theta^2}{2}.
		$$
		Moreover, since $\gamma_i^{l+1}\leq \max\{L_0,L\}$ and
		$\lambda^{l+1}\in\Delta_m$,
		$$
		\sum_{l=N_2}^{k-1}
		\sum_{i=1}^m
		\frac{\lambda_i^{l+1}}{\gamma_i^{l+1}}
		\geq
		\frac{k-N_2}{\max\{L_0,L\}}.
		$$
		Define
		$
		\bar{\mu}_i^{N_2,k}
		:=
		\frac{
			\sum_{l=N_2}^{k-1}
			\lambda_i^{l+1}/\gamma_i^{l+1}
		}{
			\sum_{j=1}^m
			\sum_{l=N_2}^{k-1}
			\lambda_j^{l+1}/\gamma_j^{l+1}
		},
		\quad i\in[m].
		$
		Then $\bar{\mu}^{N_2,k}\in\Delta_m$, and hence
\[
    \min_{i\in[m]}\{f_i(x^k)-f_i(y)\} \leq\sum_{i=1}^m\bar\mu_i^{N_2,k} \bigl(f_i(x^k)-f_i(y)\bigr) \leq\frac{\alpha\max\{L_0,L\}\Theta^2}{2(k-N_2)}.
    \]
		Taking the supremum over $y\in\Omega_k$ yields
		$$
		u_0(x^k)
		\leq
		\frac{\alpha\max\{L_0,L\}\Theta^2}{2(k-N_2)}.
		$$
	\end{proof}

	\section{Numerical experiments}

	\subsection{Some implementation details}
	We compare AC-SDM and AC-NGDM with the multiobjective steepest descent method using an Armijo line search (MSD--Armijo) and a max-type nonmonotone line search (MSD-max type) on the test problems listed in Table~\ref{1}. All experiments were conducted in MATLAB R2023b on a PC with an AMD Ryzen 5 5600 processor and 16 GB of RAM.

\begin{table}[!htbp]
    \caption{In the Convex column, $\checkmark$ denotes a convex problem and $\times$ denotes a nonconvex problem. The $L_0$ values apply to both AC-SDM and AC-NGDM. The listed dimensions and boxes are those used in this study; $x_L$ and $x_U$ generate the initial points.}
    \label{1}\label{tab:test-problems}
    \fontsize{9}{10}\selectfont
		\centering
		\renewcommand{\arraystretch}{1.0}
		\setlength{\tabcolsep}{3pt}
			\begin{tabularx}{\textwidth}{@{}lcccccc>{\raggedright\arraybackslash}X@{}}
				\toprule
				Problem & Convex & $n$ & $m$ & $L_0$ & $x_L$ & $x_U$ & Reference \\
				\midrule
				DD1 & $\times$ & 5 & 2 & $0.1$ & $[-20,\ldots,-20]$ & $[20,\ldots,20]$ & \citep{chen2023bb} \\
				DGO1 & $\times$ & 1 & 2 & $0.5$ & $[-10]$ & $[13]$ & \citep{huband2006test} \\
				FDS & $\checkmark$ & 10 & 3 & $0.5$ & $[-2,\ldots,-2]$ & $[2,\ldots,2]$ & \citep{chen2023bb} \\
				FF1 & $\times$ & 2 & 2 & $0.1$ & $[-1,-1]$ & $[1,1]$ & \citep{huband2006test} \\
				Imbalance2 & $\checkmark$ & 2 & 2 & $2$ & $[-2,-2]$ & $[2,2]$ & \citep{chen2023bb} \\
				JOS1a & $\checkmark$ & 50 & 2 & $0.1$ & $[-2,\ldots,-2]$ & $[2,\ldots,2]$ & \citep{chen2023bb} \\
				JOS1b & $\checkmark$ & 100 & 2 & $0.1$ & $[-2,\ldots,-2]$ & $[2,\ldots,2]$ & \citep{chen2023bb} \\
				JOS1c & $\checkmark$ & 100 & 2 & $0.1$ & $[-50,\ldots,-50]$ & $[50,\ldots,50]$ & \citep{chen2023bb} \\
				JOS1d & $\checkmark$ & 100 & 2 & $0.1$ & $[-100,\ldots,-100]$ & $[100,\ldots,100]$ & \citep{chen2023bb} \\
				KW2 & $\times$ & 2 & 2 & $0.1$ & $[-3,-3]$ & $[3,3]$ & \citep{mita2019} \\
				LTY3a & $\times$ & 100 & 2 & $0.1$ & $[-2,\ldots,-2]$ & $[2,\ldots,2]$ & \citep{yin2025flow} \\
				LTY3b & $\times$ & 100 & 2 & $0.1$ & $[-2,\ldots,-2]$ & $[2,\ldots,2]$ & \citep{yin2025flow} \\
				LTY3c & $\times$ & 100 & 2 & $0.1$ & $[-2,\ldots,-2]$ & $[2,\ldots,2]$ & \citep{yin2025flow} \\
				LTY3d & $\times$ & 100 & 2 & $0.1$ & $[-2,\ldots,-2]$ & $[2,\ldots,2]$ & \citep{yin2025flow} \\
				MLF1 & $\times$ & 1 & 2 & $0.5$ & $[0]$ & $[20]$ & \citep{huband2006test} \\
				MOP2 & $\times$ & 3 & 2 & $0.1$ & $[-4,-4,-4]$ & $[4,4,4]$ & \citep{huband2006test} \\
				MOP7 & $\checkmark$ & 2 & 3 & $0.1$ & $[0,0]$ & $[10,10]$ & \citep{huband2006test} \\
				PNR & $\times$ & 2 & 2 & $2$ & $[-2,-2]$ & $[2,2]$ & \citep{preuss2006pareto} \\
				QV1 & $\times$ & 50 & 2 & $0.1$ & $[-5.12,\ldots,-5.12]$ & $[5.12,\ldots,5.12]$ & \citep{huband2006test} \\
				SD & $\checkmark$ & 4 & 2 & $0.4$ & $[1,\sqrt{2},\sqrt{2},1]$ & $[3,3,3,3]$ & \citep{mita2019} \\
				TKLY1 & $\times$ & 4 & 2 & $0.1$ & $[0.1,0,0,0]$ & $[1,1,1,1]$ & \citep{huband2006test} \\
				TOI4a & $\checkmark$ & 50 & 2 & $1.5$ & $[-2,\ldots,-2]$ & $[5,\ldots,5]$ & \citep{mita2019} \\
				TOI4b & $\checkmark$ & 100 & 2 & $1.5$ & $[-2,\ldots,-2]$ & $[5,\ldots,5]$ & \citep{mita2019} \\
				TOI4c & $\checkmark$ & 200 & 2 & $1.5$ & $[-2,\ldots,-2]$ & $[5,\ldots,5]$ & \citep{mita2019} \\
				TOI4d & $\checkmark$ & 500 & 2 & $1.5$ & $[-2,\ldots,-2]$ & $[5,\ldots,5]$ & \citep{mita2019} \\

				TRIDIA2 & $\times$ & 4 & 4 & $8$ & $[-1,\ldots,-1]$ & $[1,\ldots,1]$ & \citep{chen2023bb} \\
				VNT2 & $\checkmark$ & 2 & 3 & $0.1$ & $[-4,-4]$ & $[4,4]$ & \citep{viennet1996multicriteria} \\
				WFG4 & $\times$ & 24 & 3 & $0.1$ & $[0,0,\ldots,0]$ & $[2,4,\ldots,48]$ & \citep{huband2006test} \\
				WIT4 & $\checkmark$ & 2 & 2 & $3$ & $[-2,-2]$ & $[2,2]$ & \citep{chen2023bb} \\
				WIT5 & $\checkmark$ & 2 & 2 & $2$ & $[-2,-2]$ & $[2,2]$ & \citep{chen2023bb} \\
				\bottomrule
			\end{tabularx}

	\end{table}

\FloatBarrier

	For each test problem, the four algorithms were tested from the same set of 200 distinct initial points, with a stopping tolerance of \(\varepsilon=10^{-4}\) and a maximum of 1000 iterations. We set \(\alpha=1.1\) for AC-SDM and AC-NGDM, with the corresponding values of \(L_0\) reported in Table~\ref{1}. For MSD--Armijo, we used \(\beta=0.5\) and \(\rho=10^{-4}\). For MSD-max type, we used \(\beta=0.5\), \(\rho=10^{-4}\), \(M=4\).

\subsection{Results analysis}
\begin{figure}[H]
  \centering
  \includegraphics[width=\linewidth]{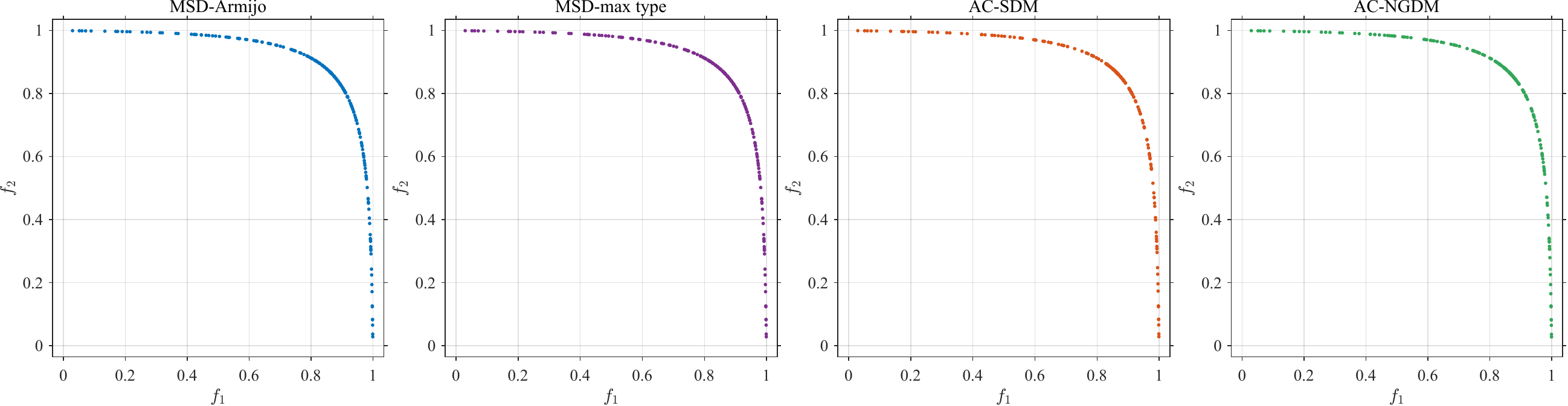}
  \caption{Numerical results in value space obtained by MSD-Armijo, MSD-max type, AC-SDM and AC-NGDM for problem FF1.}
  \label{fig:objective-FF1}
\end{figure}

\begin{figure}[H]
  \centering
  \includegraphics[width=\linewidth]{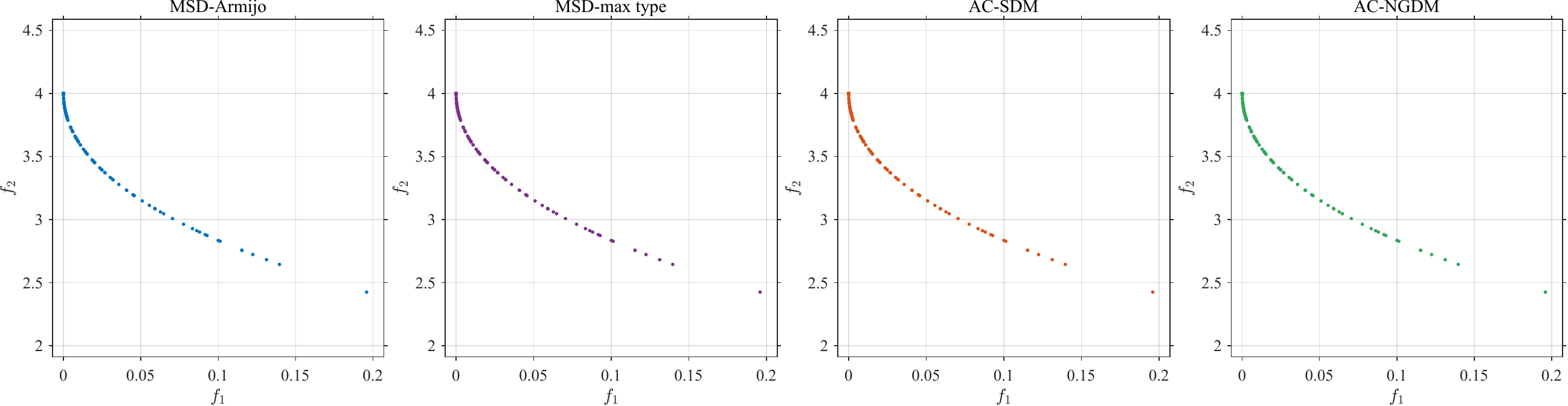}
  \caption{Numerical results in value space obtained by MSD-Armijo, MSD-max type, AC-SDM and AC-NGDM for problem JOS1a.}
  \label{fig:objective-JOS1a}
\end{figure}

\begin{figure}[H]
  \centering
  \includegraphics[width=\linewidth]{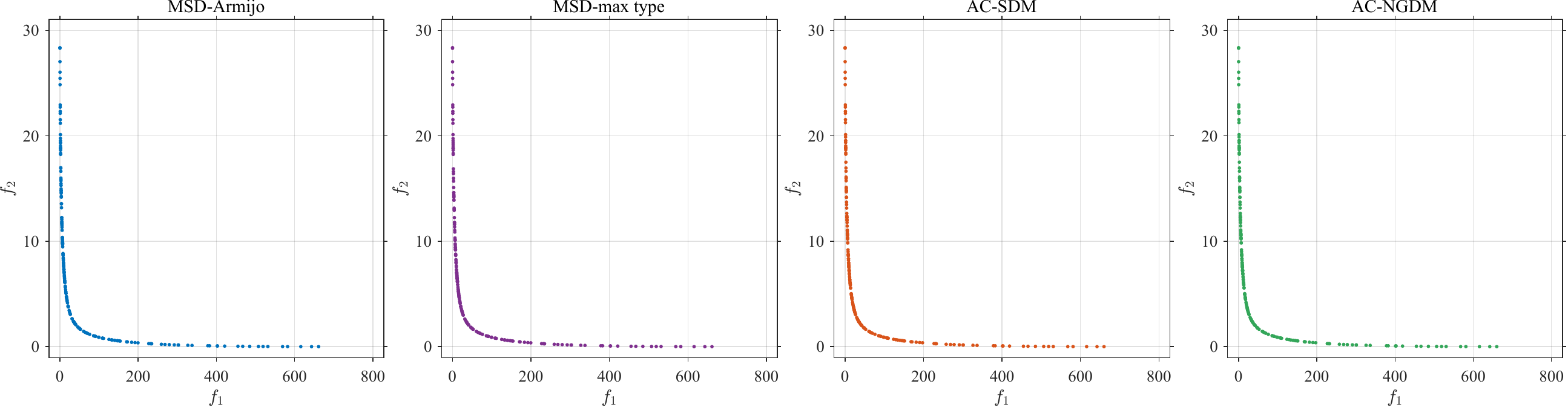}
  \caption{Numerical results in value space obtained by MSD-Armijo, MSD-max type, AC-SDM and AC-NGDM for problem WIT4.}
  \label{fig:objective-WIT4}
\end{figure}

\begin{figure}[H]
  \centering
  \includegraphics[width=\linewidth]{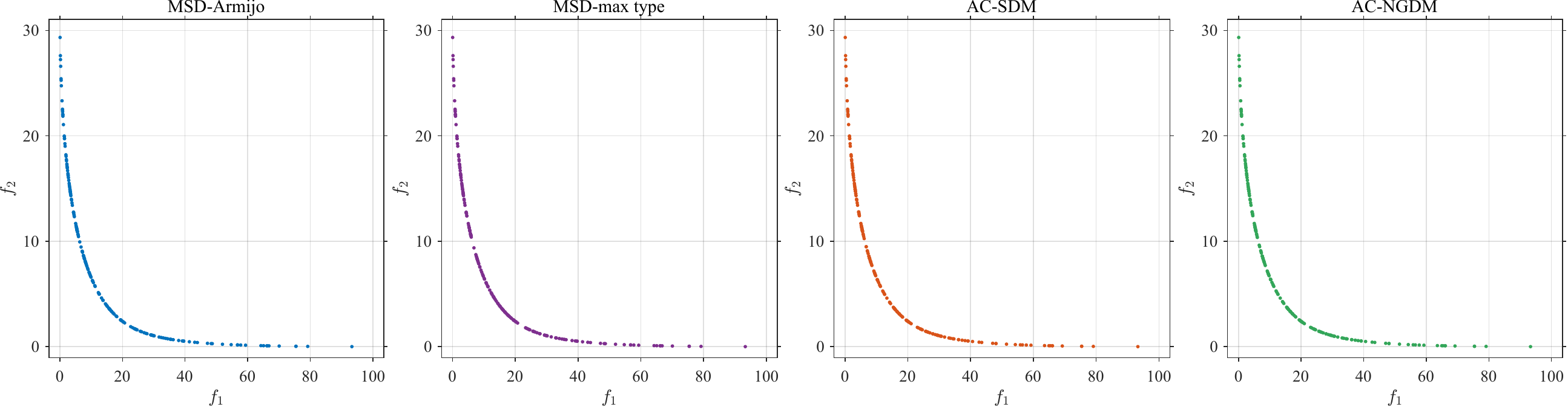}
  \caption{Numerical results in value space obtained by MSD-Armijo, MSD-max type, AC-SDM and AC-NGDM for problem WIT5.}
  \label{fig:objective-WIT5}
\end{figure}

\begin{figure}[H]
  \centering
  \includegraphics[width=\linewidth]{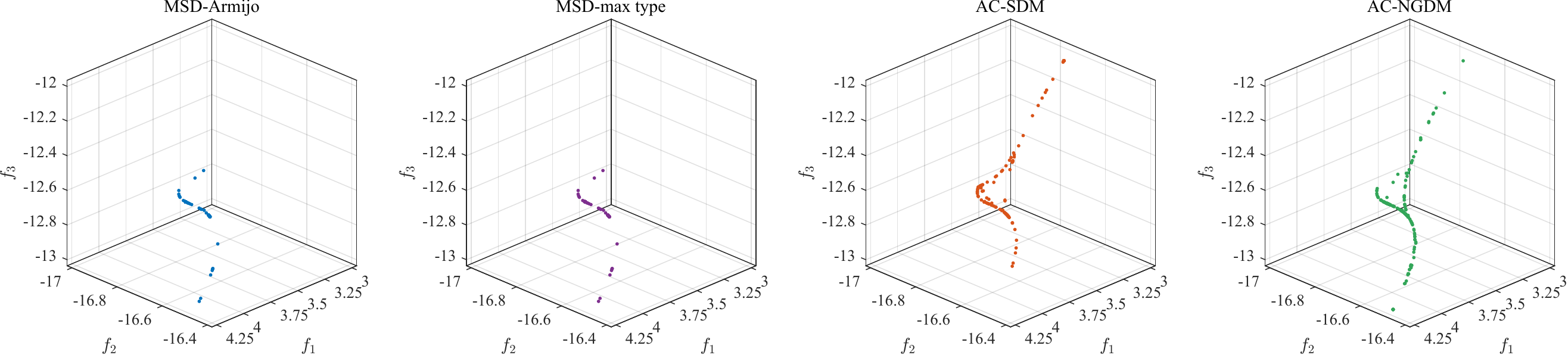}
  \caption{Numerical results in value space obtained by MSD-Armijo, MSD-max type, AC-SDM and AC-NGDM for problem MOP7.}
  \label{fig:objective-MOP7}
\end{figure}

\begin{figure}[H]
  \centering
  \includegraphics[width=\linewidth]{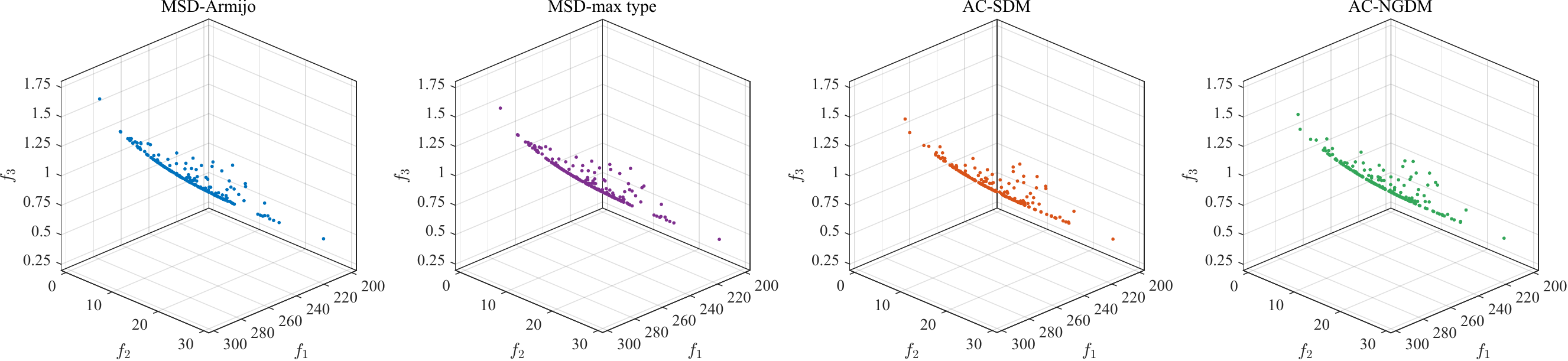}
  \caption{Numerical results in value space obtained by MSD-Armijo, MSD-max type, AC-SDM and AC-NGDM for the quadratic FDS problem.}
  \label{fig:objective-FDS}
\end{figure}

The numerical comparisons are reported in Table~\ref{tab:common-advantage-results}. Figures~\ref{fig:objective-FF1}--\ref{fig:objective-FDS}
show the corresponding value space.

\begin{table}[!htbp]
    \caption{Time(s), iter, and conv(\%) denote the total CPU time in seconds, the average iteration count over recorded runs, and the convergence rate as a percentage of all 200 initial points, respectively.}
    \label{tab:common-advantage-results}
    \fontsize{7}{8.5}\selectfont
		\centering
		\renewcommand{\arraystretch}{1.0}
		\setlength{\tabcolsep}{1pt}
			\begin{tabular}{@{}l@{\hspace{3pt}}rrr@{\hspace{3pt}}rrr@{\hspace{3pt}}rrr@{\hspace{3pt}}rrr@{}}
				\toprule
				Problem
				& \multicolumn{3}{c}{MSD-Armijo}
				& \multicolumn{3}{c}{MSD-max type}
				& \multicolumn{3}{c}{AC-SDM}
				& \multicolumn{3}{c}{AC-NGDM} \\
				\cmidrule{2-4}\cmidrule{5-7}\cmidrule{8-10}\cmidrule{11-13}
				& time(s) & iter & \shortstack{conv.\\(\%)}
				& time(s) & iter & \shortstack{conv.\\(\%)}
				& time(s) & iter & \shortstack{conv.\\(\%)}
				& time(s) & iter & \shortstack{conv.\\(\%)} \\
				\midrule
				DD1 & 4.703 & 69.00 & \textbf{100\%} & 5.641 & 36.09 & \textbf{100\%} & \textbf{3.813} & 37.26 & \textbf{100\%} & 4.484 & \textbf{9.01} & 99.5\% \\
				DGO1 & 0.453 & 2.49 & \textbf{100\%} & 0.609 & 2.49 & \textbf{100\%} & \textbf{0.156} & \textbf{1.18} & \textbf{100\%} & 0.172 & \textbf{1.18} & \textbf{100\%} \\
				FDS & 11.484 & 101.62 & \textbf{100\%} & 5.094 & 56.75 & \textbf{100\%} & 8.906 & 125.02 & \textbf{100\%} & \textbf{3.094} & \textbf{37.32} & \textbf{100\%} \\
				FF1 & 3.063 & 43.57 & \textbf{100\%} & 3.563 & 43.57 & \textbf{100\%} & 3.219 & 39.33 & \textbf{100\%} & \textbf{1.250} & \textbf{7.73} & \textbf{100\%} \\
				Imbalance2 & 24.625 & 305.55 & 98.5\% & 18.359 & 288.80 & \textbf{100\%} & 61.734 & 1000.00 & 0\% & \textbf{0.375} & \textbf{4.28} & \textbf{100\%} \\
				JOS1a & 13.000 & 198.45 & \textbf{100\%} & 13.031 & 198.45 & \textbf{100\%} & \textbf{1.563} & \textbf{18.22} & \textbf{100\%} & 1.969 & 23.47 & \textbf{100\%} \\
				JOS1b & 25.703 & 383.61 & \textbf{100\%} & 24.469 & 383.61 & \textbf{100\%} & \textbf{3.219} & \textbf{39.02} & \textbf{100\%} & 3.844 & 50.63 & \textbf{100\%} \\
				JOS1c & 34.875 & 542.98 & \textbf{100\%} & 34.250 & 542.98 & \textbf{100\%} & \textbf{3.969} & \textbf{55.03} & \textbf{100\%} & 5.000 & 66.67 & \textbf{100\%} \\
				JOS1d & 37.172 & 577.35 & \textbf{100\%} & 35.938 & 577.35 & \textbf{100\%} & \textbf{4.328} & \textbf{58.63} & \textbf{100\%} & 5.266 & 70.02 & \textbf{100\%} \\
				KW2 & 11.625 & 179.05 & 87.5\% & 14.375 & 229.21 & 87.5\% & \textbf{4.922} & \textbf{69.01} & \textbf{96.5\%} & 10.094 & 148.64 & 89.5\% \\
				LTY3a & 1.844 & \textbf{6.87} & \textbf{100\%} & 11.656 & 177.90 & \textbf{100\%} & 1.938 & 20.95 & \textbf{100\%} & \textbf{1.391} & 18.40 & \textbf{100\%} \\
				LTY3b & 1.516 & 13.04 & \textbf{100\%} & 1.563 & 13.04 & \textbf{100\%} & 1.047 & 7.15 & \textbf{100\%} & \textbf{0.750} & \textbf{6.10} & \textbf{100\%} \\
				LTY3c & 1.438 & 12.74 & \textbf{100\%} & 1.594 & 12.73 & \textbf{100\%} & 1.422 & 11.06 & \textbf{100\%} & \textbf{0.938} & \textbf{6.65} & \textbf{100\%} \\
				LTY3d & 1.422 & 13.09 & \textbf{100\%} & 1.594 & 13.06 & \textbf{100\%} & 1.656 & 13.86 & \textbf{100\%} & \textbf{1.172} & \textbf{7.77} & \textbf{100\%} \\
				MLF1 & 0.141 & 0.86 & \textbf{100\%} & 0.203 & 0.86 & \textbf{100\%} & \textbf{0.125} & \textbf{0.68} & \textbf{100\%} & 0.141 & \textbf{0.68} & \textbf{100\%} \\
				MOP2 & 2.953 & 41.44 & 98\% & 3.188 & 41.44 & 98\% & \textbf{0.859} & \textbf{6.51} & \textbf{100\%} & 3.781 & 50.24 & 97\% \\
				MOP7 & 4.641 & 69.80 & \textbf{100\%} & 4.813 & 69.80 & \textbf{100\%} & \textbf{2.828} & 33.90 & \textbf{100\%} & 3.984 & \textbf{32.85} & \textbf{100\%} \\
				PNR & 1.516 & 9.24 & \textbf{100\%} & 0.844 & 7.82 & \textbf{100\%} & 4.156 & 56.43 & \textbf{100\%} & \textbf{0.766} & \textbf{5.51} & \textbf{100\%} \\
				QV1 & 67.219 & 1000.00 & 0\% & 65.188 & 1000.00 & 0\% & \textbf{44.063} & \textbf{662.53} & \textbf{96\%} & 52.531 & 771.50 & 88\% \\
				SD & 1.734 & 22.94 & \textbf{100\%} & 1.781 & 22.94 & \textbf{100\%} & 2.063 & 16.47 & \textbf{100\%} & \textbf{1.422} & \textbf{14.33} & \textbf{100\%} \\
				TKLY1 & 34.844 & 532.99 & 60\% & 32.047 & 498.61 & 62.5\% & \textbf{20.656} & \textbf{318.14} & 75\% & 23.734 & 354.71 & \textbf{92.5\%} \\
				TOI4a & 16.156 & 237.24 & 77\% & 27.328 & 403.55 & 60.5\% & \textbf{0.797} & 6.09 & \textbf{100\%} & \textbf{0.797} & \textbf{5.67} & \textbf{100\%} \\
				TOI4b & 5.000 & 67.10 & 99.5\% & 3.391 & 44.75 & \textbf{100\%} & 0.703 & 5.53 & \textbf{100\%} & \textbf{0.672} & \textbf{5.16} & \textbf{100\%} \\
				TOI4c & 3.969 & 47.12 & \textbf{100\%} & 2.656 & 29.51 & \textbf{100\%} & \textbf{0.734} & 5.71 & \textbf{100\%} & 0.922 & \textbf{5.43} & \textbf{100\%} \\
				TOI4d & 3.906 & 43.38 & \textbf{100\%} & 2.969 & 27.39 & \textbf{100\%} & \textbf{1.031} & 5.93 & \textbf{100\%} & 1.969 & \textbf{5.83} & \textbf{100\%} \\

				TRIDIA2 & 4.813 & 43.60 & 99.5\% & 3.406 & 42.00 & 99.5\% & 6.328 & 90.41 & 99.5\% & \textbf{3.109} & \textbf{37.62} & \textbf{100\%} \\
				VNT2 & 6.219 & 89.46 & \textbf{100\%} & 6.125 & 89.46 & \textbf{100\%} & 4.125 & 52.25 & \textbf{100\%} & \textbf{2.781} & \textbf{29.35} & \textbf{100\%} \\
				WFG4 & 187.094 & 43.31 & \textbf{100\%} & 1228.719 & 317.81 & 92.5\% & 165.031 & 38.73 & \textbf{100\%} & \textbf{133.406} & \textbf{30.76} & \textbf{100\%} \\
				WIT4 & 0.922 & 8.60 & \textbf{100\%} & 1.828 & 18.39 & \textbf{100\%} & 1.141 & 15.98 & \textbf{100\%} & \textbf{0.625} & \textbf{7.66} & \textbf{100\%} \\
				WIT5 & 1.344 & 14.00 & 99.5\% & 2.516 & 32.76 & 99\% & 0.641 & 9.67 & \textbf{100\%} & \textbf{0.609} & \textbf{4.17} & \textbf{100\%} \\
				\bottomrule
			\end{tabular}

	\end{table}

For the convex problems JOS1a–JOS1d and TOI4a–TOI4d, both AC-SDM and AC-NGDM use less CPU time and take fewer iterations than MSD–Armijo and MSD-max type. With the problem dimension increasing, the advantages of the two proposed methods become more evident. For example, AC-SDM takes 4.328s in total on JOS1d and averages 58.63 iterations. The corresponding average is 577.35 for both line-search methods, with total CPU times of 37.172s for MSD–Armijo and 35.938 s for MSD-max type. Both proposed methods meet the stopping criterion in every run on TOI4a–TOI4d, taking approximately five to six iterations. MSD-max type takes fewer iterations on average than MSD–Armijo on TOI4b–TOI4d, but still more than AC-SDM and AC-NGDM. AC-NGDM has the lowest total CPU time and average iteration count on FDS, VNT2, WIT4, and WIT5. Across the four methods, the total CPU time ranges from 2.828s to 4.813s for MOP7 and from 1.422s to 2.063s for SD. In particular, Imbalance2 is a bi-objective function in which the gradient norms of the two objective functions differ substantially. The performance of AC-NGDM on this problem illustrates its advantage in handling such imbalances.

For the nonconvex problem QV1, neither line-search method meets the stopping criterion within 1000 iterations in any of the 200 runs. AC-SDM and AC-NGDM meet it in 96\% and 88\% of the runs, respectively. On TKLY1, both proposed methods also use less CPU time than line-search method. The percentage of runs satisfying the stopping criterion is highest for AC-NGDM at 92.5\%, followed by AC-SDM at 75\%, MSD-max type at 62.5\%, and MSD–Armijo at 60\%. AC-SDM requires the least CPU time and the fewest iterations on average for KW2 and MOP2, while AC-NGDM does so for FF1, PNR, and TRIDIA2. The total CPU times on PNR are close for MSD-max type and AC-NGDM, at 0.844s and 0.766s, respectively.

\FloatBarrier

\section{Conclusion}

In this paper, we proposed AC-SDM and AC-NGDM, two methods that use local estimates of the Lipschitz constants to construct descent directions without line-search. AC-SDM constructs a scaled steepest descent direction using a shared estimate, whereas AC-NGDM constructs its direction from gradients normalized by corresponding estimates. We proved that all accumulation points of the sequences generated by AC-SDM and AC-NGDM are Pareto critical points and established convergence rates for both nonconvex and convex cases. Numerical experiments indicate that both methods perform favorably compared to MSD-Armijo and MSD-max type on the test problems. 

It is notable that the current choice of \(\gamma\) involves taking a maximum and may therefore be strongly influenced by an excessively large local Lipschitz estimate for a single objective. Although AC-NGDM alleviates this issue, more flexible choices of \(\gamma\) remain worth investigating to further reduce this effect. Besides, AC-SDM and AC-NGDM can be extended to proximal gradient methods for composite multiobjective optimization problems.


\end{document}